\documentclass[review]{elsarticle}

\usepackage{amsmath}
\usepackage{amsthm}
\usepackage{graphicx}
\usepackage{fullpage}
\usepackage[utf8]{inputenc}

\newtheorem{theorem}{Theorem}

\newtheorem{conj}{Conjecture}

\DeclareMathOperator*{\argmin}{arg\,min}

\begin{document}

\newcommand{\RR}{\ensuremath{\mathbf{R}^2}}
\newcommand{\N}{\ensuremath{\matbf N}}
\newcommand{\Z}{\ensuremath{\mathbf Z}}

\begin{frontmatter}

\title{On the numerical stability of the least-squares method for the planar scattering by obstacles}
\author{Gilles Chardon \fnref{adress}}
\ead{gilles.chardon@m4x.org}
\fntext[address]{tel. +43 1 51581-2500, fax. +43 1 51581 2530}
\address{Acoustics Research Institute, Austrian Academy of Sciences, Wohllebengasse 12-14, 1040 Wien, Österreich}


\begin{abstract}
The scattering of waves by obstacles in
a 2D setting is considered, in particular the computation of the scattered field
via the collocation or the least-squares methods. In the case of multiple
scattering by
smooth obstacles, we prove that the scattered field
can be uniformly approximated by sums of multipoles. For a unique obstacle, the choice of the number of points
and their positions for the estimation of the error on the border of the scatterer is studied,
showing the benefit of using a non-uniform distribution of points
dependent on the scatterer and the approximation scheme. In general,
using a denser discretization near the singularities of the scattered
field does not improve the stability of the method.
The analysis
can also be used to estimate the discretization size needed
to ensure stability given a density of points and an approximation scheme,
e.g. in the case of multiple scatterers.
\end{abstract}

\begin{keyword}
least-squares method, scattering, numerical quadrature, Helmholtz equation, Trefftz methods
\end{keyword}

\end{frontmatter}


\section{Introduction}

This article is concerned with the least-squares method
for the Helmholtz equation \cite{Monk, Stojek}.
Like other methods such as the Boundary Element
Method, the Variational Theory of Complex Rays \cite{vtcr}, or the Ultra-weak
variational Formulation \cite{uwvf}, it solves
the Helmholtz equation by using an approximation scheme 
for the solutions to the equation (plane waves,
generalized harmonic polynomials, etc.). These methods
differ by the way they match the solution to the boundary conditions
and ensure the continuity between the elements.
With the least-squares method, the boundary conditions and the continuity between
the subdomains are enforced via the minimization of a $L_2$-norm, allowing
a simple implementation of the method.

In the particular case of the scattering of an incident wave $u_i$ by obstacles $S$ in the plane with Dirichlet boundary conditions, the scattered field is solution to:
\begin{equation}
\left\{
\begin{array}{l}
\Delta u_s + k^2 u_s = 0\\
u_s = - u_i \mbox{ on } \partial S\\
\lim_{r\rightarrow \infty} \sqrt{r} \left( \frac{\partial u_s}{\partial r} - ik u_s\right) = 0
\end{array}
\right.
\label{helm}
\end{equation}
The resolution of this problem via the least-squares method was studied, among others, by Stojek \cite{Stojek}, and Barnett and Betcke \cite{Barnett}. The domain of propagation
was partitioned in bounded subdomains where Fourier-Bessel functions,
or fractional Fourier-Bessel functions in the case of a scatterer with corners,
were used for the approximation, and an unbounded domain.
In this domain, Stojek used Hankel functions to enforce
the Sommerfeld radiation condition, while Barnett and Betcke
used the Method of Fundamental Solutions.
We restrict ourselves to the somewhat simpler
case, but equally interesting, of the scattering by smooth obstacles,
that, as will be shown in this article, allows
 the approximation of the solutions on the entire unbounded domain
of propagation using a single set of functions.

In the case of a unique scatterer, this approximation is the basis of the so-called Rayleigh methods\cite{millar}, where
the scattered field is approximated by sums of multipoles
\begin{equation}
u_s \approx u_N = \sum_{n = -N_h}^{N_h} \alpha_n H_n(kr)e^{in\theta}
\end{equation}
where $(r, \theta)$ are the polar coordinates.
 The collocation method, or Point Matching Method, estimates the $2N_h+1$ coefficients
by fitting the boundary conditions on $2N_h+1$ points on the border of the scatterer.
While simple to implement, this method is usually numerically unstable as the matrix to be inverted
is likely to be ill-conditioned.  The coefficients can also be estimated by minimizing the $L_2$ error on the
boundary, where the error is approximated by numerical quadrature.  This can be considered as a particular case of the least-squares methods cited above,
where the approximation is done on the entire domain like in \cite{eisenstat}.
However, in the case of scattering by an obstacle, the domain is multiply-connected and unbounded. As is usually observed in this particular case \cite{kleev, semenova}, but
also in the general case \cite{Stojek, Monk, Barnett}, the matrices
involved in the computation have the tendency to be ill-conditioned.

In this article, we investigate the effect of the choice of the quadrature points on the numerical stability of the least-squares
method in the particular case of the scattering by an obstacle.
The quadrature used to estimate the error is generally either
left unspecified \cite{Monk, Stojek, ramm, desmet}, or uses general purpose
schemes (Chebyshev nodes in \cite{betcke} or Clenshaw-Curtis rule in
\cite{Barnett}).
However, as was shown in \cite{kleev} for the collocation method,
the choice of the quadrature point is critical for the stability
of the numerical methods, and depends on the shape of the scatterer.
This has direct implications on the efficiency of the computational methods, as
choosing the appropriate quadrature rule allows to use fewer points, making
the matrices involved in the computation smaller.

This work is a first step towards a more general study of the effect of the quadrature scheme on the numerical stability of least-squares methods (i.e.
with other approximations schemes such as plane waves, Fourier-Bessel functions,
fractional Fourier-Bessel functions, and with several subdomains),
but is also interesting in itself as it gives a stable numerical
scheme for the scattering by an obstacle, a long standing problem
in electromagnetics and acoustics\cite{millar,semenova, kleev, christiansen}.

In section 2, we prove an approximation result for the scattering of waves by smooth obstacles:
given a set of $L$ smooth scatterers $S_l$, the scattered field $u_s$, solution
to \eqref{helm}, can be uniformly approximated
by a sequence of sums of multipoles
\begin{equation}
u^N =  \sum_{l = 1}^L\sum_{n\in \Z} \beta_{nl}^N H_n(kr_l)e^{in\theta_l}
\end{equation}
where $(r_l, \theta_l)$ are the polar coordinates with respect to a center $O_l$, with at least one such center in each scatterer, and a finite number
of coefficients $ \beta_{nl}^N$ are nonzero.
In the  third section of the paper, the stability of the least-squares method is investigated in the light
of a recent result by Cohen et al. \cite{cdl} on least-squares approximations, and numerical results
are given for the scattering by an ellipse and a square, demonstrating the importance of the choice
of the quadrature points. In particular, it is shown that in some
cases, using a denser
discretization near the singularities makes the stability harder to achieve
than using a uniform density or even a denser discretization away from the
singularities. Finally, the application to the scattering of
a plane wave by two arbitrary shaped scatterers is given.

\section{Approximation of the scattered field}

We here show that the field scattered by smooth obstacles can be uniformly approximated by sums of multipoles. This is an application of the Vekua theory \cite{vekua, henricivekua}, a theory of elliptic
PDEs allowing the construction of operators mapping holomorphic functions to solutions of a given PDE, as long
as the coefficients of the PDE are analytic. This is obviously true for the Helmholtz equation, for
which the operator mapping holomorphic (or equivalently harmonic) functions to solutions to the Helmholtz equation
and its inverse are explicitly known. These operators being continuous, approximation results
available for holomorphic functions can be translated to similar results for solutions to the Helmholtz equation. For instance, in a star-shaped domain,
solutions to the Helmholtz equation can be approximated by sums of Fourier-Bessel functions $J_n(kr) e^{in\theta}$, that are the images of the harmonic polynomials used
to approximate harmonic functions.
More details on the Vekua operators
and on the approximation of solutions to the Helmholtz equation in convex domains can be found 
in the articles by Moiola et al. \cite{moiolavekua, moiola2}.

\begin{figure}
\centering
\includegraphics[width=7cm]{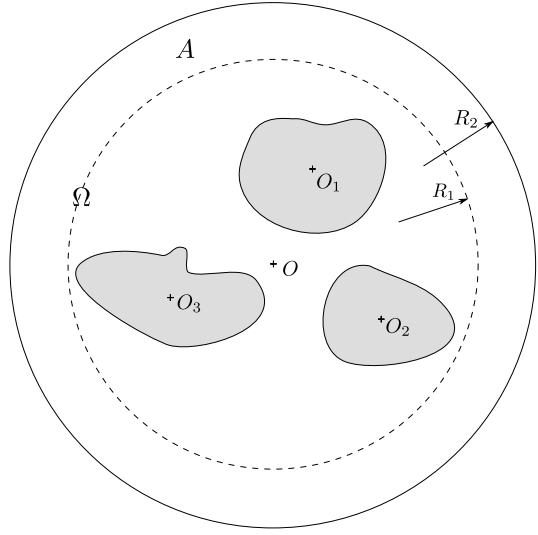}
\caption{Domain on which the Vekua theory is applied}
\label{domain}
\end{figure}

Our setting is as follows: scatterers are contained in the disk of radius
$R_1$, and $\Omega$ is the closed domain delimited by the circle of radius $R_2 > R_1$
and the scatterers, see figure~\ref{domain}.

\begin{theorem}
Let $S$ be a set of scatterers in the plane, and $u_s$ the scattered field.
If $u_s$ can be analytically continued in a open domain containing
$\mathbf{R}^2-S$, then $u_s$ can be uniformly approximated, in $\mathbf{R}^2-S$
as well as on the boundaries of $S$, by a sequence $u^N$
of sums of Fourier-Hankel functions 
\begin{equation}
u^N =  \sum_{l = 1}^L\sum_{n\in \Z} \beta_{nl}^N H_n(kr_l)e^{in\theta_l}
\end{equation}
where $(r_l, \theta_l)$ are the polar coordinates associated to the centers $O_l$,
with at least such a center in each scatterer, and finitely many
coefficients  $\beta_{nl}^N$ are nonzero for a given $N$.
\label{theoapp}
\end{theorem}

\begin{proof}
The proof is divided in three main parts:
\begin{itemize}
\item the approximation of the scattered field in $\Omega$ by sums of first-kind Fourier-Bessel functions and outgoing Fourier-Hankel functions,
\item the use of the Sommerfeld radiation condition to reduce
this approximation to sums of outgoing Hankel functions,
\item the extension of this approximation outside of $\Omega$.
\end{itemize}

\paragraph{General representation in $\Omega$}

For a closed multiply-connected domain $\Omega$ and the Helmholtz equation, Vekua proved (see Ref.\cite{vekua} page 109) that if a solution $u$ can be analytically continued in an open domain containing $\Omega$,
then $u$  can be uniformly approximated by a sequence of finite sums of Fourier-Bessel functions of the first and second kind
$$J_n(kr)e^{in\theta}, Y_n(k_l r)e^{in\theta_l}$$
where $(r, \theta)$ are the polar coordinates with respect to the origin $O$, and $(r_l, \theta_l)$
with respect to $O_l$, with at least a point $O_l$ arbitrarily chosen in the $l$-th 
simply connected component of $\mathbf{R}^2-\Omega$. As we can uniformly approximate the functions $J_n(kr_l)e^{in\theta_l}$ by sums of $J_m(kr)e^{im\theta}$ (Graf theorem), we can equivalently approximate $u$ by sums of first kind Fourier-Bessel functions and outgoing Fourier-Hankel functions:
$$J_n(kr)e^{in\theta}, H_n(k_l r)e^{in\theta_l}$$
We apply this theorem in $\Omega$. We have thus that $u^N$, defined
as 
\begin{equation}
u^N
 = \sum_{n\in \Z}\alpha_n^N J_n(kr) e^{in\theta} + \sum_{l = 1}^L\sum_{n\in \Z} \beta_{nl}^N H_n(kr_l)e^{in\theta_l}
\label{vekua}
\end{equation}
where only a finite number of coefficients $\alpha_n^N$ and $\beta_{nl}^N$ are nonzero,
uniformly converges to $u_s$ as $N \rightarrow \infty$.

\paragraph{Removal of the first term in (\ref{vekua})}

In the annulus $A$ of radiuses $R_1$ and $R_2$  we can, by moving
the Hankel functions to $O$ (using the Graf theorem), approximate $u_s$ as the limit of a sum of first-kind Bessel functions and a series of Hankel functions
with coefficients functions of $N$.
\begin{equation}
u^N = \sum_{n\in \Z} \alpha_n^{N} J_n(kr)e^{in\theta} + \beta_n^N H_n(kr)e^{in\theta}
\label{seq}
\end{equation}
For $r \geq R_1$, $u$ can also be written as a series of Hankel functions
as it satisfies the Sommerfeld radiation conditions:
\begin{equation}
u_s = \sum_{n\in \Z} \gamma_n H_n(kr) e^{in\theta}
\label{series}
\end{equation}
Fitting the sequence \eqref{seq} to the series \eqref{series} in the annulus $A$ will allow to show that $u_s$
can be approximated by Hankel functions around the scatterers.

Let $\epsilon > 0$ and $R'_2$ with $R_1 < R'_2 < R_2$. Then there is a $N_1$ such that for $N > N_1$,
\begin{equation}
|u^N - u_s| = \left|\sum_{n\in \Z} \alpha_n^{N} J_n(kR'_2)e^{in\theta} + \beta_n^N H_n(kR'_2)e^{in\theta}
- \gamma_n H_n(kR'_2)e^{in\theta}\right| < \epsilon
\end{equation}
The coefficients of the Fourier series of $u^N - u_s$  satisfy
\begin{equation}
 \left| \alpha^{N}_n J_n(kR'_2)  + \beta^N_n H_n(kR'_2)- \gamma_n H_n(kR'_2)\right| < \epsilon
\end{equation}
As the functions $u^N$ and their limit are analytic, we can
do the same for the radial derivative. For $N > N_2$:
\begin{equation}
 k\left| \alpha^{N}_n J_n'(kR'_2) +  \beta^N_n H_n'(kR'_2) - \gamma_n H_n'(kR'_2)\right| < k \epsilon
\end{equation}

For $N > \max(N_1, N_2)$, we have
\begin{equation}
\alpha^{N}_n J_n(kR'_2) + \beta^N_n H_n(kR'_2) = \gamma_n H_n(kR'_2) + \delta_n
\end{equation}
\begin{equation}
\alpha^{N}_n J_n'(kR'_2) + \beta^N_n H_n'(kR'_2) = \gamma_n H_n'(kR'_2) + \delta'_n
\end{equation}
with $|\delta_n| < \epsilon$ and $|\delta'_n| < \epsilon$

Solving this system for $\alpha^{N}_n$  (using the 
fact that the wronskian\cite{AS} of $J_n$ and $H_n$ is equal to $2i/(\pi kR'_2$)
we have
\begin{equation}
 \alpha^{N}_n = \frac{2 i(\delta_n H_n'(kR'_2) - \delta_n' H_n(kR'_2))}{\pi kR'_2},
\end{equation}
and, with $C = 2/\pi kR'_2$
\begin{equation}
 |\alpha^{N}_n| \leq  C\epsilon \left(|H_n'(kR'_2)| + |H_n(kR'_2)|\right).
\end{equation}

We now uniformly bound the first term of \eqref{seq} for $r \leq R''_2 <  R'_2 $. We have
$|\sum \alpha^{N}_n J_n(kr)| \leq \sum |\alpha^{N}_n J_n(kr)| \leq C \epsilon \sum 
(  |H_n'(kR'_2)|+ |H_n(kR'_2)|) |J_n(kr)|$ and
\begin{equation}
\sum_{n \in \Z} (|H_n'(kR'_2)|+ |H_n(kR'_2)|)| J_n(kr)| 
= \sum_{|n| < kR''_2} (|H_n'(kR'_2)|+ |H_n(kR'_2)|)| J_n(kr)| + \sum_{|n| \geq kR''_2} (|H_n'(kR'_2)|+ |H_n(kR'_2)|)| J_n(kR''_2)| 
\label{truc}
\end{equation}
where we used the fact that $|J_n|$ is increasing on $[0, n]$.
This first term is a continuous function and can be uniformly bounded, and the second term is a convergent series:
\begin{align}
(|H_n'(kR'_2)|+ |H_n(kR'_2)|)|J_n(kR''_2)| & \sim |H_{n+1}(kR'_2)||J_n(kR''_2)|\\
& \sim \frac{1}{\sqrt{n+1}} \left(\frac{ekR'_2}{2(n+1)}\right)^{-n-1}
\frac{1}{\sqrt{n}} \left(\frac{ekR''_2}{2n}\right)^{n}\\
& \sim (R''_2/R'_2)^n / (kR'_2)
\end{align}
We thus have, with $C'$ the quantity in \eqref{truc},
\begin{equation}
\left|\sum_{n\in\Z} \alpha^{N}_n J_n(kr)\right| \leq \epsilon CC'
\end{equation}
The first term of \eqref{vekua} can be made uniformly as small as desired,
and can be removed from the sequence $u^N$ without changing its limit
for $ r < R''_2$.

\paragraph{Convergence outside $\Omega$}

We now have that the sequence
\begin{equation}
u^N =  \sum_{l = 1}^L\sum_{n\in \Z} \beta_{nm}^{N} H_n(kr_l)e^{in\theta_l}
\end{equation}
 uniformly converges to $u_s$ outside of
the obstacles and on their boundaries, for $r \leq R''_2$.

For $R_1 \leq r \leq R_2$, we have (cf. \eqref{seq}) that
\begin{equation}
u^N = \sum_{n\in \Z} \beta_n^N H_n(kr)e^{in\theta}
\end{equation}
converges uniformly to
\begin{equation}
u_s = \sum_{n\in \Z} \gamma_n H_n(kr) e^{in\theta}.
\end{equation}

 On the circle of radius
$R_1$, the error between $u_N$ and $u$ converges uniformly to 0.
The coefficients of the Fourier transform of the error can be bounded:
\begin{equation}
|\beta_n^N - \gamma_n| < \epsilon / H_n(kR_1)
\end{equation}
Now, for any $r > R'_1$ (with $R_1 < R'_1 < R''_2$), we have
\begin{align}
| u^N - u_s| &\leq \sum_{n \in \Z} |\beta_n^N - \gamma_n| |H_n(kr)\|\\
& \leq \epsilon \sum_{n \in \Z} |H_n(kR'_1)| / |H_n(kR_1)|\\
& \leq \epsilon C
\end{align}
as the general term of the last series is equivalent to $(R'_1/R_1)^n$.
We here used the fact that $|H_n|$ is a decreasing function.

The error between $u^N$ and $u_s$ can be bounded by any $\epsilon$
for any $r > R'_1$. Combined with the fact that $u^N$ converges
uniformly to $u_s$ for $r < R''_2$ outside of the obstacles
and on their boundaries,
$u^N$ converges uniformly to $u_s$ outside of the scatterers
and on their boundaries, yielding the theorem.
\end{proof}

To conclude this section, we formulate a conjecture based on the analogy between Runge's theorem and Theorem \ref{theoapp}. 
This theorem states that an analytic function on a given open domain can be,
in a compact subdomain, uniformly approximated by a sequence of
rational functions, with at least a pole in each connected component
of the complement of its domain of analyticity. The connection
between this result for analytic functions and solutions 
to the Helmholtz equation is given by the Vekua theory.
 The conditions of  Runge's theorem can
actually be weakened, as shown by Mergelyan's theorem\cite{greene}. This theorem states that it is sufficient
that the function to be approximated is holomorphic in the interior
of $\Omega$ and only continuous on the boundaries of $\Omega$. In particular, singularities can be present on the boundaries, which
is likely for domains with corners. It is then
reasonable to formulate this conjecture:

\begin{conj}
The uniform approximation by sum of multipoles is valid as long as the scattered field is analytic in the exterior
of the obstacles, and continuous on their boundaries.
\end{conj}

This condition is in particular satisfied for domain with corners, Dirichlet
conditions and a continuous incident field. Although this
conjecture is of theoretical interest, its usefulness in numerical applications is limited.
Indeed, the presence of singularities on the boundary of the obstacles
limits the rate of convergence of the Vekua approximations \cite{moiola2}. Convergence
can be accelerated by the use of fractional Fourier-Bessel functions in the approximation\cite{Barnett},
but this necessitates the partition of the exterior domain in simply connected subdomains, making multipole approximations irrelevant.

\section{Stability of numerical methods}

In the previous section, we proved that uniform approximation of the scattered
field was possible for smooth scatterers, that is that the error between
the scattered field $u_s$ and its best approximation $u^N$ of order $N$ tends 
uniformly to zero:
$$\| u_s - u^N\|_\infty \rightarrow 0.$$

We are here interested in least-squares method for the case of a unique
scatterer. The uniform convergence
implies local convergence in the $L_2$-norm (that
is, on any compact domain):
$$\| u_s - u^N\|_{2,loc} \rightarrow 0.$$
In particular, it implies the convergence in the $L_2$ norm on the
boundary $\Gamma$ of the scatterer. The following theorem, proved
in \cite{ramm}, shows that the convergence on the boundary of the scatterer
is sufficient to ensure convergence outside:

\begin{theorem}
(Ramm, Gutman) If $g$ is in $L^2(\Gamma)$, then the solution $w$
of the Helmholtz equation with Sommerfeld radiation conditions
and $w = g$ on $\Gamma$ is bounded on the exterior domain $D'$ by
\begin{equation}
\|w\| \leq C \| g\|_{L^2(\Gamma)}
\end{equation}
where $\|w\| = \|w\|_{H^m_{loc}(D')} +   \|w\|_{L^2(D', (1+|x|)^{-\gamma})}$,
with $H^m$ with $m>0$ is the Sobolev space and the $L^2$ norm is weighted
by $(1+|x|)^{-\gamma}$ with $\gamma >1$.
\label{ramm}
\end{theorem}

However, this does not guarantee that practical
computation of the scattered field will always converge to the 
true solution when increasing the order of approximation.
To this end,  the estimation of the coefficients of this finite approximation
has to be  stable, that is that the error (e.g. on the boundary $\Gamma$) between this estimated field $\tilde u^N$ and the
actual field is of the same order as the best approximation error:
$$\| u_s - \tilde u^N\|_{2,\Gamma} \approx \| u_s - u^N\|_{2,\Gamma}.$$

In this section, we analyze the numerical stability of the computation
of the scattered field using the multipole approximation and 
the collocation or least-squares methods. Through numerical evaluation of the stability, we aim to show that the density of points used on the border
of the scatterers is critical for the stability. The tool
used for the analysis can also be used to estimate, given
an approximation scheme and a sampling density (e.g. multipole
approximation with uniform density on the border), the number
of samples needed to ensure stability.

With the least-squares methods, the scattered field is estimated as follows: given an order of approximation $N_h$ (i.e. $m = 2N_h +1$ Fourier-Hankel functions) and
a number $N_s$ of points on the boundary of the scatterer, the coefficients
of the multipole approximation are estimated by matching the incoming field $u_i$ and the multipole approximation on the sampling points in the least-squares sense:
\begin{equation}
\tilde{\mathbf a} = \argmin_{\mathbf a} \| \mathbf u + \mathbf H \mathbf a \|_2
\end{equation}
where $\tilde{\mathbf a} = (\alpha_{-N_h}, \ldots , \alpha_{N_h})$ is the vector containing the estimated coefficients, $\mathbf u$
the vector of the values of $u_{i}$ sampled on the boundary, and $\mathbf H$ the $N_s \times 2N_h+1$ matrix
with terms
$$H_{mn} = H_n(kr_m) e^{in\theta_m}$$
where $n \in \{-N_h, \ldots N_h\}$ and $(r_m, \theta_m)$ are the polar coordinates of the $N_s$ points
on the boundary. Note that the Sommerfeld radiation condition
does not need to be considered in the minimization problem
as it is enforced through the choice of the approximation spaces.
The estimated scattered field is then given by
$$\tilde u_{s} =  \sum_{|n| \leq N_h} \alpha_n H_n(kr) e^{in\theta}.$$

The collocation method, or point matching method, is a particular case of the least-squares
method, obtained when $N_s = 2N_h+1$. In this case, the matrix $\mathbf H$ is square,
and the coefficients
are found by matching the values of the samples in an exact way.

With these methods, estimation of the scattered field is essentially reduced
to the interpolation of the incident field on the boundary of the
scatterer using  a finite number of functions (the traces of the multipoles on the boundary) from a finite number of samples. However, it is well known that 
even in basic cases, interpolation of a function
from a finite number of samples can be unstable (cf. Runge phenomenon), even when the data is  perfectly known on the sampling points.

\subsection{Stability of least-squares estimation}

To analyze the stability of these methods
in function of $N_h$, $N_s$ and the density of samples,
we use results
by Cohen et al \cite{cdl}. These results allow the estimation of the number
of measurements ensuring stability, knowing the density probability
measure used to draw the points on the border, as well as the desired
order of approximation. While the sampling points are not generally
chosen randomly, and the values obtained are somewhat pessimistic, these results will help to evaluate the stability
of sampling densities.

The setting is as follows. The functions to be estimated, defined on a space $X$, are known
to be approximated by elements of spaces $V_m$ of dimension $m$,
and the best approximation error of $u$ by an element of
$V_m$ is denoted by $\sigma_m(u)$. The estimation $\tilde u$ is
obtained by the least-squares method using $n$ samples,
drawn from the space $X$ using the probability density $\nu$, and
truncated so that its absolute value is not larger than
$M = \max_{x\in X} |u(x)|$.

To evaluate the stability of this least-square estimation we compute the quantity
\begin{equation}
K(m) = \max_{x \in X} \sum_{j = 1}^m | L_j(x)|^2
\end{equation}
where $(L_j)_j$ is a basis of the space $V_m$, orthogonal with respect to the probability density $\nu$:
\begin{equation}
\int_X L_j^\star L_k d\nu = \delta_{j,k}
\label{int}
\end{equation}

The value of $K(m)$ is then linked to the number of measurements $n$
by the following theorem:
\begin{theorem}
(Cohen, Davenport, Leviatan)
Let $r>0$ be arbitrary but fixed and let $\kappa:=\frac {1-\log 2}{2+2r}$. 
 If $m$ is such that 
$$K(m) \leq \kappa \frac{n}{\log n}$$
 then, the expectation of the reconstruction error is bounded:
 $$
 E(\|u-\tilde u\|^2)\leq (1+\epsilon(n))\sigma_m(u)^2+8M^2n^{-r},
 $$
 where $\epsilon(n):=\frac {4\kappa} {\log n} \to 0$ as $n\to +\infty$.
\label{theols}
\end{theorem}

In our case, the space $X$ is the border of the scatterer, and the spaces $V_m$ are generated by finite families
of multipoles. The orthogonal basis $(L_j)_j$ can
be estimated by orthogonalizing a family of multipoles, using
 the Gram-Schmidt algorithm and numerical quadrature (e.g. Monte-Carlo
integration using the probability density $\nu$).

Numerical results are now given for the scattering by an ellipse, a square,
and two ovals. The code needed to reproduce the figures, as well as
similar results for others scatterers are available online \cite{web}.

\subsection{Scattering by an ellipse}

\begin{figure}
\centering
\includegraphics[width=7cm]{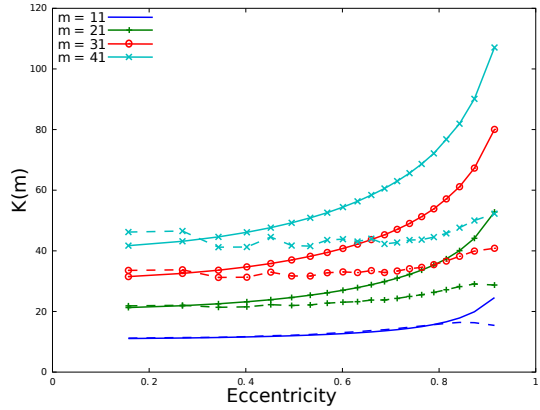}
\caption{Estimated $K(m)$ for the ellipse at $k=6$ in functions of the eccentricity, for the uniform density (solid) and the density obtained by Kleev (dashed)}
\label{fig8}
\end{figure}

We first estimate $K(m)$ for the case of the ellipse centered at the origin.
The values of $K(m)$ for $m \in \{10, 20, 30, 40, 50\}$ and $k = 6$ as a function of the eccentricity $e$ of the ellipse are given on figure \ref{fig8}, for the uniform density and the density suggested by Kleev and Manenkov\cite{kleev} .
This density (denoted by KM density in the rest of the paper) is based on a conformal mapping between the unit disk and
the scatterer. The computation of this density is outlined in the appendix.

\begin{figure}
\centering
\includegraphics[width=7cm]{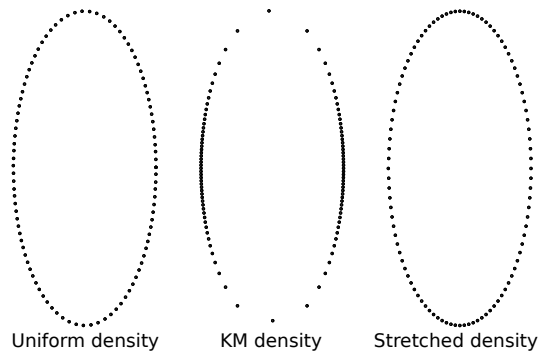}
\caption{Samplings on the ellipse, with eccentricity $e = 0.95$. Left: uniform density, center: KM density, right: uniform density on the circle mapped to the ellipse}
\label{fig_ellipse}
\end{figure}

For the uniform density,
the value of $K(m)$ increases with the eccentricity, meaning
that more and more points are needed for a fixed number of coefficients.
$K(m)$ remains nearly constant for the KM density.
 It is clear that this density needs less samples to ensure the stability
of the interpolation. 
However, in contrast to the claim of Kleev and Manenkov that the density depends on the singularities of the scattered field,
these results shows that the density does not depend on these singularities, but on the singularities of the
functions used for the interpolation. Indeed, as the eccentricity of the
ellipse increases, the singularities, located at the focal points
of the ellipse, approaches the extremities of the major axis, while
the sampling becomes denser near the extremities of the minor axis,
i.e. near the singularities of the functions used to approximate the scattered field (see figure \ref{fig_ellipse}).

\begin{figure}
\centering
\includegraphics[width=7cm]{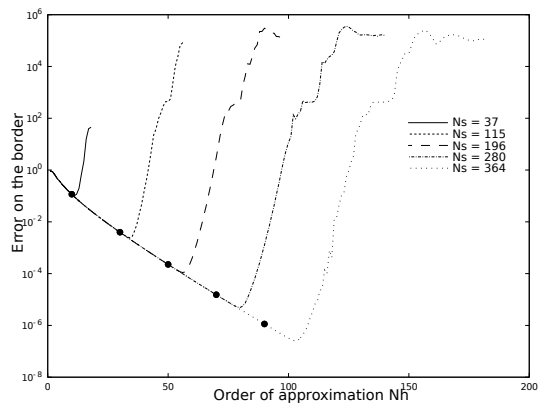}
\caption{$L_2$-Error on the boundary for the scattering of a plane wave by an ellipse,
with eccentricity $e = 0.8$ for several numbers of samples, and varying orders
of approximation. For each number of samples, the order $N_h$ such that
$N_s \approx K(2N_h+1)$ is highlighted.}
\label{fig_K_ellipse}
\end{figure}

On figure \ref{fig_K_ellipse}, we plot, for some fixed numbers $N_s$ of samples
on the boundary, the $L_2$ approximation error on the border in function
of the approximation order $N_h$ for $e = 0.8$ and $k = 6$.
This error on the border is an indicator of the quality of the estimation
of the scattered field. Indeed, as Theorem \ref{ramm} shows, the error outside of the scatterer can be bounded
by the error on the border.
The orders for which $N_s = K(2N_h+1)$ are indicated. As the approximation
errors for these choices of parameters are close the optimal
errors, we suggest to use $K(2N_h+1)$ samples when using an order
of approximation of $N_h$.

\begin{figure}
\centering
\includegraphics[width=7cm]{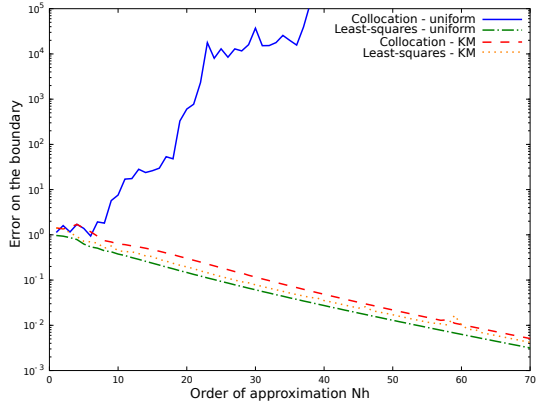}
\caption{$L_2$-error on the boundary for the scattering of a plane wave by an ellipse
with eccentricities $e = 0.86$ and varying approximation orders. The error is given for the collocation and the
least-squares methods, with uniform and KM densities. For the least-squares,
the number of samples is given by $K(m)$.}
\label{fig_err_ell}
\end{figure}

On figure \ref{fig_err_ell},
 the $L_2$ error between the estimated scattered field and the incident field
on the border is plotted
for the result of the collocation and the least-squares method, with the uniform and the KM densities,
in function of the approximation order for an eccentricity $e = 0.86$ and $k = 10$.
For the least-squares estimation, we use $K(m)$ samples.
The collocation with
uniform density fails as the error increases with the approximation order.
Using the KM density makes the collocation method stable. For the uniform density,
using $K(m)$ samples yields a stable estimation. For the KM density, using
$K(M)$ samples slightly improves the results.

We now estimate $K(m)$ when the scattered field is approximated
using Mathieu functions, which give separable solutions
to the Helmholtz equation in elliptic coordinates \cite{barakat}. 
This method is used in \cite{lee} to compute the scattering
by multiple ellipses. We test
here two densities, the uniform density, and the density obtained
by stretching the uniform density on a circle (see Fig. \ref{fig_ellipse}). Note that the Mathieu
functions are orthogonal for this second density, and that their values
on the ellipse do not depend on the wavenumber. We find here that in this case, a better stability is obtained
by using the stretched density, that is when using more samples near the large axis. This is coherent with the observations above as in this case, the singularities of the functions (products of Mathieu functions in elliptical coordinates)
are at the focal points of the ellipse, near the end of the major axis.

\begin{figure}
\centering
\includegraphics[width=7cm]{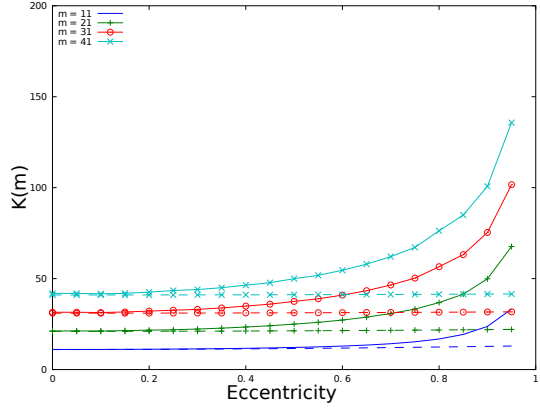}
\caption{Estimated $K(m)$ for the ellipse in functions of the eccentricity, using Mathieu functions, for the uniform density (solid) and the uniform density on the circle mapped to the ellipse (dashed)}
\label{fig4}
\end{figure}

\subsection{Scattering by a square}

\begin{figure}
\centering
\includegraphics[width=7cm]{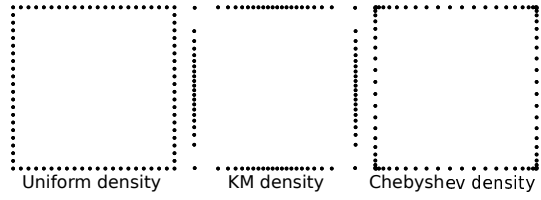}
\caption{Samplings on the square}
\label{fig7}
\end{figure}

\begin{figure}
\centering
\includegraphics[width=7cm]{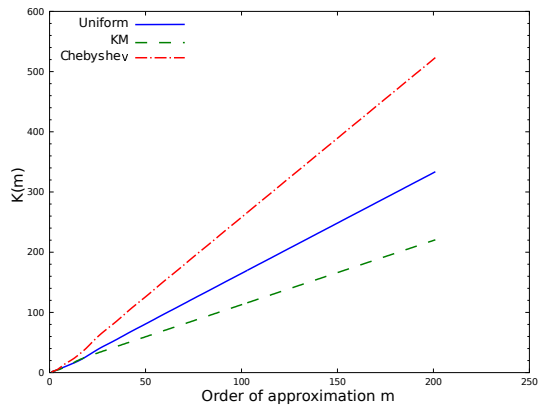}
\caption{Estimated $K(m)$ for the square at $k=6$ for the uniform,
KM and Chebyshev densities.}
\label{ksq}
\end{figure}

We now give results for the case of the unit square. We test here three densities (see figure \ref{fig7}):
\begin{itemize}
\item the uniform density,
\item the KM density,
\item the density $1/4\pi\sqrt{1-\min(x, y)^2}$.
\end{itemize}
The last density (called Chebyshev density in the rest of the paper),
is similar to the sampling given by the Chebyshev nodes, used in the Clenshaw-Curtis quadrature rule.
The estimation of $K(m)$ for these three densities is plotted
on figure \ref{ksq}. Although the singularities of the scattered
field are localized at the corners of the square, using
a denser discretization near those corners is actually harmful
to the stability of the least-squares method. Surprisingly, using
more samples near the center of the edges of the square, i.e.
far from the singularities, yields a $K(m)$ slightly larger
than $m$, its lower bound.

\begin{figure}
\centering
\includegraphics[width=7cm]{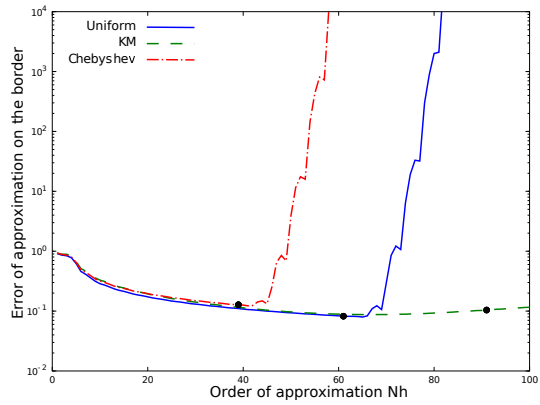}
\caption{$L_2$-Error on the boundary for the scattering of a plane wave by a square for the three densities on the border and varying orders
of approximation, $k = 6$. For each density, the order $N_h$ such that
$N_s \approx K(2N_h+1)$ is highlighted.}
\label{ksqerr}
\end{figure}

The $L_2$ error on the border for a fixed number of samples and the
three densities is plotted on figure \ref{ksqerr} as a function
of the order of approximation. Like above, the order
for which $K(2N_h +1) \approx N_s$ is indicated. The KM density is stable
for any approximation order, in particular for the collocation case. The uniform density can yield
a slightly lower error, but is not always stable. Using
more points near the corners does not allow to use
a large order of approximation, and gives the largest error.

\subsection{Scattering by two ovals}

\begin{figure}
\centering
\includegraphics[width=7cm]{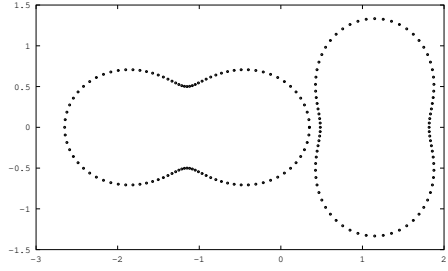}
\caption{Samplings on the ovals}
\label{fig10}
\end{figure}

Our final numerical experiments deals with the scattering by two Booth ovals,
defined in polar coordinates by
$$r(\theta) = 1 + \cos(2\theta)/a.$$
As indicated by Theorem 1, the scattered field can be approximated by
two families of multipoles on the entire domain of propagation. The centers of the multipoles are chosen at the centers of the ovals (with parameters $a = 2$ and $a = 3$), and their borders
are sampled by using a uniform sampling of the angle $\theta$, plotted on figure
\ref{fig10}. To ensure the stability of the least-squares method, we
compute $K(m)$ for the two ovals separately. With $N_h = 65$, i.e. $m = 131$, we
find $K(m) = 328$ for $a = 2$ and $K(m) = 233$ for $a = 3$.

\begin{figure}
\centering
\includegraphics[width=12cm]{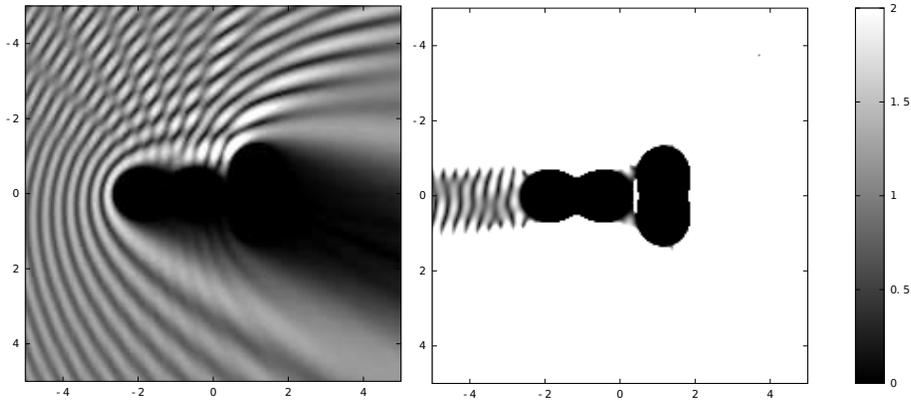}
\caption{Scattering of a plane wave by two ovals (absolute value).
Left: least-squares method. Right: collocation (thresholded at 2).}
\label{fig11}
\end{figure}

The scattered field for a incident plane wave (incoming at an angle 0.3 from the x-axis) is pictured on figure \ref{fig11}.
The total number of degrees of freedom is 262, and 561 samples are used.
Although the values of $K(m)$ are computed for the two ovals separately,
using these values in the case of multiple scattering yields a stable
estimation of the scattered field. This is expected, as the instabilities
are mostly caused by the high order Hankel functions, which are decaying
rapidly. The influence of such a function associated to a scattered
on the other scatterer is thus negligible. For comparison, the result
of the collocation method (i.e. using 131 samples on each scatterer) is also given.

\section{Conclusion}

The computation of the acoustical field scattered by obstacles was considered, in particular the numerical stability of the least-squares method.
We proved that the field scattered by smooth obstacles can be uniformly
approximated by sums of multipoles, and that the choice of the centers of the multipoles is only limited by the constraint that at least 
one center is chosen in each scatterer.

We also investigated the stability of the least-squares 
method based on
multipole approximations.
This stability depends on
 the set of functions used
to approximate the scattered field, and on the density of samples
used on the boundary of the scatterer. In particular, it does
not depends on the location of the singularities of the functions
to be approximated. Using more points near these singularities
can even be detrimental to the stability.
We showed that a simple numerical computation can yield,
given a set of functions and a density of samples, an estimate
of the number of samples necessary to ensure stability of the least-squares method, and that it can be also used in the case of multiple scatterers.
In the case of multipole approximations, the densities suggested by Kleev and Manenkov are close to the optimum.

On a more general level, we showed that in the context of
least-squares or collocation methods, the approximation scheme
and the quadrature rule have to be chosen conjointly.

\section{Acknowledgments}

The author is supported by the Austrian Science Fund (FWF) START-project FLAME (Frames and Linear Operators for Acoustical Modeling and Parameter Estimation; Y 551-N13), and thanks Vincent Pagneux and Alexandre Leblanc for fruitful discussions.


\appendix
\section{Determination of the KM points}

The densities suggested by Kleev and Manenkov are obtained
by a mapping the exterior of the unit disk to the image of
the scatterer by inversion. The KM points
are the images of a uniform sampling of the disk by this mapping.
Note that this is equivalent to a mapping from the interior of the
unit disk to the scatterer with the origin as fixed point.

\subsection{Case of the ellipse}

A conformal mapping from the unit disk to the ellipse of major semi-axis $a$
and minor semi-axis 1 is given by (see \cite{ellipse})
$$f(z) =  \sqrt{a^2-1} \sin \left( \frac{\pi}{2 K(k)}
\int_0^{z/\sqrt{k}} \frac{dt}{\sqrt{(1-t^2)(1-k^2t^2)}}\right)$$
where $K$ is the complete elliptic integral of the first kind. The parameter $k$
is found as the solution of
$$
\frac{K(k')}{K(k)} = \frac{2}{\pi} \mathrm{asinh}\left(\frac{2a}{a^2-1}\right)
\mbox{ and } k' = \sqrt{1-k^2}.$$

The KM points are then simply
the images of $e^{i2\pi n/N}$ by $f$.

\subsection{Case of the square}

A conformal mapping from the unit disk to the square is given
by the Schwartz-Christoffel mapping \cite{SC} $f$, with
$$f'(z) = c \prod_{n = 1}^4  \left(1 - \frac{z}{z_k}\right)^{\alpha_k -1}$$
where $z_k$ are the inverse images of the vertices and $\pi \alpha_k$
the angles of the square. For symmetry reason, we choose $z_1 = 1$,
$z_2 = i$, $z_2 = -1$, $z_2 = -i$. The $\alpha_k$ are equal to $1/2$. The map
$$f(z) = \int_0^z \frac{1}{\sqrt{1-z'^4}} dz'$$
thus maps the unit disk to a square.

\end{document}